\documentclass[a4paper,12pt]{article}
\usepackage{amssymb}

\def\coloana(#1,#2){\left(\matrix{#1\cr#2\cr}\right)}
\def\matrice(#1,#2,#3,#4){\left(\matrix{#1&#2\cr#3&#4\cr}\right)}

\newtheorem{lemma}{Lemma}
\newtheorem{theorem}{Theorem}

\begin{document}

\author{Octavian G. Mustafa\\
\small{Faculty of Mathematics and Computer Science, DAL,}\\
\small{University of Craiova, Romania}\\
\small{e-mail: octawian@yahoo.com}\\
and\\
Cemil Tun\c{c}\\
\small{Department of Mathematics, Faculty of Arts and Sciences}\\
\small{Y\"{u}z\"{u}nc\"{u} Yil University, 65080 Van, Turkey}\\
\small{e-mail: cemtunc@yahoo.com}}
\title{Asymptotically linear solutions of differential equations via Lyapunov functions}
\date{}
\maketitle

\noindent\textbf{Abstract} We discuss the existence of solutions with oblique asymptotes to a class of second order nonlinear ordinary differential equations by means of Lyapunov functions. The approach is new in this field and allows for simpler proofs of general results regarding Emden-Fowler like equations.


\section{Introduction}

Let us consider the following nonlinear ordinary differential equation
\begin{eqnarray}
x^{\prime\prime}+f\left(t,\frac{x}{t}\right)=0,\qquad t\geq t_0\geq1,\label{main_eq}
\end{eqnarray}
where the nonlinearity $f:[t_0,+\infty)\times\mathbb{R}\rightarrow\mathbb{R}$ is assumed continuous. Further smoothness for the function $f$ will be added later. A model for this equation is provided by the Emden-Fowler like equation
\begin{eqnarray}
x^{\prime\prime}+A(t)x^{2n-1}=0,\qquad t\geq t_0,\label{EF_eq}
\end{eqnarray}
where $n\geq1$ is an integer.

By an \textit{asymptotically linear solution} of equation (\ref{main_eq}) we understand any $C^2$ real-valued function $x(t)$ that verifies the equation in a neighborhood of $+\infty$ and can be represented as
\begin{eqnarray}
x(t)=x_{1}t+x_{2}+o(1),\qquad x_1,x_2\in\mathbb{R},\label{AL_1}
\end{eqnarray}
together with its derivative
\begin{eqnarray}
x^{\prime}(t)=x_1+o(t^{-1})\qquad\mbox{when }t\rightarrow+\infty.\label{AL_2}
\end{eqnarray}
From now on, any $C^1$ function which is defined in a neighborhood of $+\infty$ and can be developed as in (\ref{AL_1}), (\ref{AL_2}) will be referred to as \textit{asymptotically linear}.

This type of behavior was investigated recently mostly in connection with the existence of positive solutions to a family of reaction-diffusion equations in exterior domains of $\mathbb{R}^{m}$, $m\geq2$, see \cite{Constantin1997,Agarwal_et_al_0,Agarwal_et_al,Hesaaraki,Mingarelli,Philos2007,Musta,Agarwal_et_al_1}. Details about the long time behavior of solutions to the Emden-Fowler equation (\ref{EF_eq}) can be found in the comprehensive monograph \cite{KiguradzeChanturia}.

The investigation of asymptotically linear solutions is performed in the literature of the last decade usually by means of fixed point theory. In this respect, several results regarding the equation (\ref{main_eq}) can be read in \cite{MustaRogo2002,MustaRogo}. A few studies dealt with the problem of asymptotically linear solutions by introducing a (Brauer-Trench-Wong type of) \textit{comparison technique} that leans upon the existence of a bounded maximal solution to a nonlinear first order ordinary differential equation, see the presentation in \cite[pp. 355, 357]{Agarwal_et_al}.

Our aim in this note is to complement the discussion from \cite{Agarwal_et_al_0,Agarwal_et_al} with an analysis of the asymptotically linear solutions of equation (\ref{main_eq}) by means of certain Lyapunov functions. To the best of our knowledge, this approach has not been tried before in the literature though the problem of nonoscillatory solutions has been studied via Lyapunov functions in several fundamental contributions, see \cite[Theorem 2]{Atkinson} or \cite[Theorem 2]{Wong1968}. 

\section{Preliminaries}
Let us start the computations by noticing that the equation (\ref{main_eq}) can be recast as a first order nonlinear differential system
\begin{eqnarray}
u^{\prime}=-tf(t,v),\qquad v^{\prime}=\frac{u}{t^2},\label{main_syst}
\end{eqnarray}
where
\begin{eqnarray}
\coloana(u,v)=\matrice(t,-1,0,\frac{1}{t})\coloana(x^{\prime},x),\qquad t\geq t_0.\label{sch_var}
\end{eqnarray}

We notice that the Jacobian of system (\ref{sch_var}) is $1$ which means that working with the pair $(u,v)$ will be the same as working with the classical $(x,x^{\prime})$.

The next lemma is needed in the sequel.

\begin{lemma}\label{lem1} Assume that the function $x:[t_0,+\infty)\rightarrow\mathbb{R}$ is $C^1$ and let the function $u$ be given by (\ref{sch_var}).

(i) If $u$ is bounded then $x(t)=x_{1}t+o(t)$ and $x^{\prime}(t)=x_{1}+o(1)$ when $t\rightarrow+\infty$ for a certain $x_{1}\in\mathbb{R}$.

(ii) If $\lim\limits_{t\rightarrow+\infty}u(t)=u_{\infty}\in\mathbb{R}$, then $x(t)=x_{1}t+x_{2}+o(1)$ and $x^{\prime}(t)=x_1+o(t^{-1})$ when $t\rightarrow+\infty$ for a certain $x_{2}\in\mathbb{R}$.

(iii) If $x$ is asymptotically linear then $u$ has a finite limit at $+\infty$.
\end{lemma}

\textbf{Proof} (i) Introduce an $U>0$ such that
\begin{eqnarray}
\vert u(t)\vert\leq U,\qquad t\geq t_0.\label{u_bounded}
\end{eqnarray}

We have
\begin{eqnarray*}
\vert v(T)-v(t)\vert\leq\int_{t}^{T}\vert v^{\prime}(s)\vert ds=\int_{t}^{T}\frac{\vert u(s)\vert}{s^2}ds\leq\frac{U}{t},
\end{eqnarray*}
where $T\geq t\geq t_0$. This yields the existence of
\begin{eqnarray*}
\lim\limits_{t\rightarrow+\infty}v(t)=\lim\limits_{t\rightarrow+\infty}\frac{x(t)}{t}=x_{1}\in\mathbb{R}.
\end{eqnarray*}

Following from (\ref{u_bounded}), the double inequality
\begin{eqnarray*}
\frac{x(t)-U}{t}\leq x^{\prime}(t)\leq\frac{x(t)+U}{t},\qquad t\geq t_{0},
\end{eqnarray*}
leads to $\lim\limits_{t\rightarrow+\infty}x^{\prime}(t)=x_{1}$.

(ii) We also have
\begin{eqnarray}
x(t)&=&tv(t)=t\left[v(t_0)+\int_{t_0}^{t}\frac{u(s)}{s^2}ds\right]\label{prelung_1}\\
&=&t\left[x_{1}-\int_{t}^{+\infty}\frac{u(s)}{s^2}ds\right],\nonumber
\end{eqnarray}
where $x_1=v(t_0)+\int_{t_0}^{+\infty}\frac{u(s)}{s^2}ds$. Since L'H\^{o}pital's rule implies
\begin{eqnarray*}
\lim\limits_{t\rightarrow+\infty}t\int_{t}^{+\infty}\frac{u(s)}{s^2}ds=\lim\limits_{t\rightarrow+\infty}u(t)=u_{\infty},
\end{eqnarray*}
we conclude that
\begin{eqnarray*}
x(t)=t\left[x_{1}-(u_{\infty}+o(1))\cdot\frac{1}{t}\right]=x_{1}t+x_{2}+o(1)\qquad\mbox{when }t\rightarrow+\infty,
\end{eqnarray*}
where $x_{2}=-u_{\infty}$. 

Further,
\begin{eqnarray*}
x^{\prime}(t)&=&x_1-\int_{t}^{+\infty}\frac{u(s)}{s^2}ds+\frac{u(t)}{t}\\
&=&x_1+t^{-1}\cdot \xi(t),\nonumber
\end{eqnarray*}
where $\xi(t)=u(t)-t\int_t^{+\infty}\frac{u(s)}{s^2}ds$. Given an $\varepsilon>0$ there exists a $T_{\varepsilon}>t_0$ such that $\vert u(t)-u(s)\vert\leq\varepsilon$ for any $s\geq t\geq T_{\varepsilon}$. This yields
\begin{eqnarray*}
\vert \xi(t)\vert\leq t\int_t^{+\infty}\frac{\vert u(t)-u(s)\vert}{s^2}ds\leq\varepsilon,\qquad t\geq T_{\varepsilon}.
\end{eqnarray*}

(iii) From the relations (\ref{AL_1}), (\ref{AL_2}) we get
\begin{eqnarray*}
u(t)&=&tx^{\prime}(t)-x(t)\\
&=&t[x_1+o(t^{-1})]-[x_{1}t+x_{2}+o(1)]\\
&=&-x_{2}+o(1)\qquad\mbox{when }t\rightarrow+\infty,
\end{eqnarray*}
that is, $\lim\limits_{t\rightarrow+\infty}u(t)=-x_{2}\in\mathbb{R}$.

The proof is complete. $\square$

Let $x$ be a solution of the equation (\ref{main_eq}) with the maximal interval of existence $[t_0,T_{\infty})$ for some $T_{\infty}\leq+\infty$ and assume that the associated function $u$ from (\ref{sch_var}) is bounded on $[t_0,T_{\infty})$. 

Then, the relation (\ref{prelung_1}) leads to, for any $T\in(t_0,T_{\infty})$,
\begin{eqnarray*}
\vert x(t)\vert+\vert x^{\prime}(t)\vert\leq(1+T)\vert v(t_0)\vert+\frac{T+2}{t_0}\cdot\Vert u\Vert_{\infty},\quad t\in[t_0,T],
\end{eqnarray*}
and respectively to
\begin{eqnarray*}
\limsup\limits_{t\nearrow T_{\infty}}[\vert x(t)\vert+\vert x^{\prime}(t)\vert]<+\infty,
\end{eqnarray*}
which means that the solution $x$ does not ``explode'' in finite time. Wintner's non-local existence theorem implies that this solution $x$ of equation (\ref{main_eq}) is defined throughout $[t_0,+\infty)$.

The heart of our approach consists of two steps. First, we use Lyapunov like functions to establish that there exist solutions to (\ref{main_syst}) which are \textit{Lagrange stable} (bounded near $+\infty$). Lemma \ref{lem1} (i) shows that it is enough to establish the boundedness of $u$ in this respect. Second, using a sign condition for $f$, we deduce the eventual monotonicity of $u$. Lemma \ref{lem1} (ii) will imply then that these (Lagrange stable) solutions correspond to asymptotically linear solutions of (\ref{main_eq}).

\section{Asymptotic integration of equation (\ref{main_eq})}

The main results of this note are the following theorems, applicable to the case of equation (\ref{EF_eq}).
\begin{theorem}\label{th1}
Suppose that the nonlinearity $f(t,v)$ of the equation (\ref{main_eq}) satisfies the conditions below
\begin{eqnarray}
vf(t,v)\leq0,\qquad\frac{\partial f}{\partial v}(t,v)\leq0,\label{hyp_f_1}
\end{eqnarray}
the function $\frac{\partial f}{\partial v}(t,v)$ being continuous in $[t_0,+\infty)\times\mathbb{R}$, and
\begin{eqnarray}
\vert f(t,v)\vert\leq a(t)g(\vert v\vert),\qquad t\geq t_0,v\in\mathbb{R},\label{hyp_f_2}
\end{eqnarray}
where the functions $a:[t_0,+\infty)\rightarrow[0,+\infty)$, $g:[0,+\infty)\rightarrow[0,+\infty)$ are continuous, $g(\alpha)>0$ for any $\alpha>0$, $g$ is monotone nondecreasing and
\begin{eqnarray}
K=\int_{t_0}^{+\infty}ta(t)dt<+\infty,\qquad\int_{1}^{+\infty}\frac{d\xi}{g(\xi)}<+\infty.\label{main_hyp}
\end{eqnarray}
Then, for any solution $x$ of equation (\ref{main_eq}) such that
\begin{eqnarray}
\frac{4}{t_0}\left(K+\frac{c}{g(1)}\right)<\int_{1+\vert v_0\vert}^{+\infty}\frac{d\xi}{g(\xi)},\qquad c=1+\frac{u_{0}^{2}}{2},\label{hyp1}
\end{eqnarray}
where $u_{0}=u(t_0)$, $v_{0}=v(t_0)$, either $x(t)=o(t)$ or $x(t)=x_{1}t+x_{2}+o(1)$, with $x_{1}\neq0$, when $t\rightarrow+\infty$.
\end{theorem}

\textbf{Proof.} Introduce the function
\begin{eqnarray*}
V(t,u,v)=\frac{u^2}{2}+u\int_{t_0}^{t}sf(s,v)ds,\qquad t\geq t_0,u,v\in\mathbb{R}.
\end{eqnarray*}

Given $(u,v)$ a local solution of the differential system (\ref{main_syst}), we compute the total derivative of $V$ with respect to $t$, that is
\begin{eqnarray*}
&&\frac{d}{dt}[V(t,u(t),v(t))]=\frac{\partial V}{\partial t}+\frac{\partial V}{\partial u}u^{\prime}+\frac{\partial V}{\partial v}v^{\prime}\\
&&=u(t)tf(t,v(t))+\left[u(t)+\int_{t_0}^{t}sf(s,v(t))ds\right]\cdot[-tf(t,v(t))]\\
&&+\left[u(t)\int_{t_0}^{t}s\frac{\partial f}{\partial v}(s,v(t))ds\right]\cdot\frac{u(t)}{t^2}\\
&&=-tf(t,v(t))\int_{t_0}^{t}sf(s,v(t))ds+\frac{[u(t)]^2}{t^2}\int_{t_0}^{t}s\frac{\partial f}{\partial v}(s,v(t))ds.
\end{eqnarray*}

Since $f(t,v)f(s,v)\geq0$ and $\frac{\partial f}{\partial v}(s,v)\leq0$ for any numbers $t\geq s\geq t_0$ and $v\in\mathbb{R}$, we conclude that
\begin{eqnarray}
\frac{d}{dt}[V(t,u(t),v(t))]\leq0\label{ineq_Lyap}
\end{eqnarray}
for as long as the solution exists to the right of $t_0$.

We integrate (\ref{ineq_Lyap}) in order to make the next estimates
\begin{eqnarray}
\frac{\vert u(t)\vert^2}{2}-y(t)\vert u(t)\vert&\leq&\frac{[u(t)]^2}{2}+u(t)\int_{t_0}^{t}sf(s,v(t))ds\nonumber\\
&\leq& V_0=V(t_0,u(t_0),v(t_0))<\vert V_0\vert+1=c,\label{quadr_ineq}
\end{eqnarray}
where
\begin{eqnarray*}
y(t)=\int_{t_0}^{t}sa(s)ds\cdot g(\vert v(t)\vert).
\end{eqnarray*}

The quadratic inequality from (\ref{quadr_ineq}), that is
\begin{eqnarray*}
\vert u(t)\vert^2-2y(t)\vert u(t)\vert-2c\leq0,
\end{eqnarray*}
 yields
\begin{eqnarray}
\vert u(t)\vert&\leq&\vert y(t)+\sqrt{[y(t)]^2+2c}\vert+\vert y(t)-\sqrt{[y(t)]^2+2c}\vert\nonumber\\
&\leq&4y(t)+2\sqrt{2c}<4[y(t)+c].\label{step1}
\end{eqnarray}

We also have, by taking into account (\ref{sch_var}),
\begin{eqnarray*}
y(t)\leq\int_{t_0}^{+\infty}\tau a(\tau)d\tau\cdot g\left(\vert v_0\vert+\int_{t_0}^{t}\frac{\vert u(s)\vert}{s^2}ds\right),\qquad v_0=v(t_0),
\end{eqnarray*}
and
\begin{eqnarray}
y(t)\leq K\cdot g(z(t)),\label{step2}
\end{eqnarray}
where
\begin{eqnarray*}
z(t)=1+\vert v_0\vert+\int_{t_0}^{t}\frac{\vert u(s)\vert}{s^2}ds.
\end{eqnarray*}

By combining (\ref{step1}), (\ref{step2}), we obtain
\begin{eqnarray*}
z^{\prime}(t)\leq\frac{4}{t^2}[Kg(z(t))+c]
\end{eqnarray*}
and
\begin{eqnarray*}
\int_{z(t_0)}^{z(t)}\frac{d\xi}{g(\xi)}&=&\int_{t_0}^{t}\frac{z^{\prime}(s)}{g(z(s))}ds\leq 4K\int_{t_0}^{t}\frac{ds}{s^2}+4c\int_{t_0}^{t}\frac{ds}{s^{2}g(z(s))}\\
&\leq&\frac{4}{t_0}\left(K+\frac{c}{g(1)}\right)<+\infty.
\end{eqnarray*}

The latter estimate implies, via the second of hypotheses (\ref{hyp1}), the boundedness of $z(t)$, $y(t)$ and $u(t)$. 

According to Lemma \ref{lem1} (i), all the solutions $x$ of equation (\ref{main_eq}) verify the asymptotic formula
\begin{eqnarray*}
x(t)=x_{1}t+o(t)\qquad\mbox{when }t\rightarrow+\infty,\qquad x_{1}\in\mathbb{R}.
\end{eqnarray*}

If $x_{1}\neq0$ then, since $\lim\limits_{t\rightarrow+\infty}v(t)=x_{1}$, the sign of $f(t,v(t))$ does not change close to $+\infty$. Consequently, the function $u$ is (non-strictly) monotonic in a neighborhood of $+\infty$. Since it is bounded, it has a finit limit at $+\infty$. The conclusion now follows from Lemma \ref{lem1} (ii). $\square$

\begin{theorem}\label{th2}
Suppose that $\frac{\partial f}{\partial t}(t,v)$ is continuous in $[t_0,+\infty)\times\mathbb{R}$. Assume also that
\begin{eqnarray*}
v f(t,v)\geq0,\qquad v\left[3f(t,v)+t\frac{\partial f}{\partial t}(t,v)\right]\leq0
\end{eqnarray*}
for all real numbers $t\geq t_0$ and $v\neq0$.

Then, if $x$ is a solution of equation (\ref{main_eq}), we have either
\begin{eqnarray*}
x(t)=o(t)\quad\mbox{or}\quad x(t)=x_{1}t+x_{2}+o(1),\;x_{1}\neq0,\quad\mbox{when }t\rightarrow+\infty.
\end{eqnarray*}
\end{theorem}

\textbf{Proof.} Introduce the function
\begin{eqnarray*}
V(t,u,v)=\frac{u^2}{2}+t^{3}\int_{0}^{v}f(t,s)ds,\qquad t\geq t_0,u,v\in\mathbb{R}.
\end{eqnarray*}

Again, the total derivative of $V$ for a local solution $(u,v)$ of the differential system (\ref{main_syst}) reads as
\begin{eqnarray*}
\frac{d}{dt}[V(t,u(t),v(t))]=t^2\int_{0}^{v(t)}\left[3f(t,s)+t\frac{\partial f}{\partial t}(t,s)\right]ds\leq0
\end{eqnarray*}
for as long as the solution $x$ of equation (\ref{main_eq}) exists to the right of $t_0$.

By an integration of this inequality, we deduce that
\begin{eqnarray*}
0\leq\frac{[u(t)]^2}{2}\leq V(t,u(t),v(t))\leq V_{0}=V(t_0,u(t_0),v(t_0))<+\infty.
\end{eqnarray*}
This estimate implies that
\begin{eqnarray*}
\vert u(t)\vert\leq U=\sqrt{2V_{0}}.
\end{eqnarray*}

The proof is completed by using Lemma \ref{lem1} and the monotonicity argument from the final part of the proof of Theorem \ref{th1}. $\square$

\section{The case of equation (\ref{EF_eq})}

\begin{theorem}
Assume that $A(t)\leq0$ for all $t\geq t_0$ and
\begin{eqnarray*}
\int_{t_0}^{+\infty}s^{2n}\vert A(s)\vert ds<+\infty.
\end{eqnarray*}
Then, for any solution $x$ of (\ref{EF_eq}) such that
\begin{eqnarray}
\frac{4}{t_0}\left(1+\frac{u_{0}^{2}}{2}+\int_{t_0}^{+\infty}t^{2n}\vert A(t)\vert dt\right)\leq\frac{1}{2n}\cdot\frac{1}{(1+\vert v_0\vert)^{2n}},\label{tehn_cond_ref_1}
\end{eqnarray}
where
\begin{eqnarray*}
u_{0}=t_{0}x^{\prime}(t_0)-x(t_0),\qquad v_{0}=\frac{x(t_0)}{t_0},
\end{eqnarray*}
either $x(t)=o(t)$ or $x(t)=x_{1}t+x_{2}+o(1)$, with $x_{1}\neq0$, when $t\rightarrow+\infty$.
\end{theorem}

\textbf{Proof.} We apply Theorem \ref{th1}. Here, $f(t,v)=t^{2n-1}A(t)\cdot v^{2n-1}$, $a(t)=t^{2n-1}\vert A(t)\vert$ and $g(\xi)={\xi}^{2n-1}$. $\square$

\begin{theorem}\label{th3}
Assume that $A(t)\geq0$ for all $t\geq t_0$ and
\begin{eqnarray}
(2n+2)A(t)+tA^{\prime}(t)\leq0,\qquad t\geq t_0.\label{Caligo}
\end{eqnarray}

Then, if $x$ is a solution of equation (\ref{EF_eq}), we have either
\begin{eqnarray}
x(t)=o(t)\quad\mbox{or}\quad x(t)=x_{1}t+x_{2}+o(1),\;x_{1}\neq0,\quad\mbox{when }t\rightarrow+\infty.\label{oblique}
\end{eqnarray}
\end{theorem}

\textbf{Proof.} We apply Theorem \ref{th2}. $\square$

\section{Comments}

The equation (\ref{EF_eq}), see \cite{MustaRogo2002}, exhibits both solutions blowing up in finite time and solutions which are not asymptotically linear. For
\begin{eqnarray*}
-A(t)=\left\{
\begin{array}{ll}
4n(2n-2)^{-2},\;t\in[1,2],\\
2n(2n-2)^{-2}(4-t),\;t\in[2,4],\\
0,\;t\geq4,
\end{array}
\right.
\end{eqnarray*}
there exists the solution $x(t)=(2-t)^{-1/(n-1)}$ defined in $[1,2)$ and for $A(t)=-2t^{2-4n}$ the equation has the solution $x(t)=t^2$, $t\geq t_0$. 

These examples prompt that it is natural to impose restrictions on the quantities $t_0$, $u_0$, $v_0$, so, conditions like (\ref{hyp1}), (\ref{tehn_cond_ref_1}) are unavoidable.

The first of conditions (\ref{main_hyp}) is necessary in the investigation of asymptotically linear solutions. In this respect, see the argumentation from \cite[Theorem III]{MooreNeh} and \cite[Section 5.2]{MustaRogo2002}. In other words, we cannot go further than the condition (\ref{main_hyp}) when looking for such a particular asymptotic behavior and this brings us to our reason for introducing the Lyapunov like quantity $V$ --- that is, to get simpler proofs of the general results of asymptotic integration for equation (\ref{main_eq}). In this respect, a result similar to our Theorem \ref{th1} that relies on a delicate application of fixed point theory can be worked out from the presentation in \cite{MustaRogo2002}.

In the case of $A(t)\geq0$, our Theorem \ref{th3} cannot cover the most flexible condition in the literature, see \cite{Waltman},
\begin{eqnarray}
\int_{t_0}^{+\infty}s^{2n}A(s)ds<+\infty.\label{Walt}
\end{eqnarray}

It provides nevertheless an interesting complement, with very simple proof, to the fundamental contribution of Atkinson \cite[Theorem 2]{Atkinson} which says that if the coefficient $A(t)>0$ verifies the Potter hypothesis (see \cite{Potter}), that is
\begin{eqnarray}
 A^{\prime}(t)\leq0,\qquad t\geq t_{0},\label{Potter}
\end{eqnarray}
together with the integral restriction
\begin{eqnarray}
\int_{t_0}^{+\infty}s^{2n-1}A(s)ds<+\infty\label{star}
\end{eqnarray}
then none of the solutions to (\ref{EF_eq}) oscillates.

The hypothesis (\ref{Caligo}) implies that 
\begin{eqnarray*}
0\leq A(t)\leq ct^{-(2n+2)},\quad t\geq t_0,
\end{eqnarray*}
for some positive number $c$, so, the conditions (\ref{Walt}), (\ref{Potter}), (\ref{star}) are verified.

According to Theorem \ref{th3}, all the solutions of equation (\ref{EF_eq}) are asymptotically linear and the solutions with $x(t)=o(t)$ when $t\rightarrow+\infty$ are nonoscillatory. The latter type of solutions exists always, see \cite[Proposition 1]{Musta}.
 
In the linear subcase ($n=1$), Theorem \ref{th3} is a variant of a technical result due to Caligo, see \cite{Caligo}.

\section{Acknowledgments}

The authors are indebted to an anonymous referee for several most useful suggestions.


\begin{thebibliography}{50}

\bibitem{Agarwal_et_al_0} R.P. Agarwal, S. Djebali, T. Moussaoui, O.G. Mustafa, Yu.V. Rogovchenko, On the asymptotic behavior of solutions to nonlinear ordinary differential equations, Asympt. Anal. 54 (2007), 1--50

\bibitem{Agarwal_et_al} R.P. Agarwal, S. Djebali, T. Moussaoui, O.G. Mustafa, On the asymptotic integration of nonlinear differential equations, J. Comput. Appl. Math. 202 (2007), 352--376

\bibitem{Agarwal_et_al_1} R.P. Agarwal, O.G. Mustafa, L. Popescu, On the positive solutions of certain semi-linear elliptic equations, Bull. Belg. Math. Soc. Simon Stevin 16 (2009), 49--57

\bibitem{Atkinson}F.V. Atkinson, On second order nonlinear oscillation, Pacific. J. Math. 5 (1955), 643--647

\bibitem{Caligo} D. Caligo, Comportamento asintotico degli integrali dell'equazione $y^{\prime\prime}(x)+A(x)y(x)=0$, nell'ipotesi $\lim_{x\rightarrow+\infty}A(x)=0$, Boll. U.M.I. 3 (1941), 286--295

\bibitem{Constantin1997}  A. Constantin, Positive solutions of quasilinear elliptic equations, J. Math. Anal. Appl. 213 (1997), 334--339

\bibitem{Hesaaraki} M. Hesaaraki, A. Moradifam, On the existence of bounded positive solutions of Schr\"{o}dinger equations in two-dimensional exterior domains, Appl. Math. Lett. 20 (2007), 1227--1231

\bibitem{KiguradzeChanturia} I.T. Kiguradze, T.A. Chanturia, Asymptotic properties of solutions of nonautonomous ordinary differential equations, Kluwer, Dordrecht, 1993

\bibitem{Mingarelli} A.B. Mingarelli, K. Sadarangani, Asymptotic solutions of forced nonlinear second order differential equations and their extensions, Electr. J. Differential Equations 2007 (2007), 1--40

\bibitem{MooreNeh} R.A. Moore, Z. Nehari, Nonoscillation theorems for a class of nonlinear differential equations, Trans. Amer. Math. Soc. 93 (1959), 30--52

\bibitem{MustaRogo2002} O.G. Mustafa, Yu.V. Rogovchenko, Global existence of solutions with prescribed asymptotic behavior for second-order nonlinear differential equations, Nonlinear Anal. TMA 51 (2002), 339--368

\bibitem{MustaRogo} O.G. Mustafa, Yu.V. Rogovchenko, Asymptotic integration of a class of nonlinear differential equations, Appl. Math. Lett. 19 (2006), 849--853

\bibitem{Musta} O.G. Mustafa, Existence of positive evanescent solutions to some quasilinear elliptic equations, Bull. Austral. Math. Soc. 78 (2008), 157--162

\bibitem{Philos2007} Ch.G. Philos, I.K. Purnaras, P.Ch. Tsamatos, Global solutions approaching lines at infinity to second order nonlinear delay differential equations, Funkc. Ekvac. 50 (2007), 213--259

\bibitem{Potter} R.L. Potter, On self-adjoint differential equations of second order, Pacific J. Math. 3 (1953), 467--491

\bibitem{Waltman} P. Waltman, On the asymptotic behavior of solutions of a nonlinear equation, Proc. Amer. Math. Soc. 15 (1964), 918--923

\bibitem{Wong1968} J.S.W. Wong, On second order nonlinear oscillation, Funkc. Ekvac. 11 (1968), 207--234

\end{thebibliography}
\end{document}